\documentclass[12pt]{amsart}
\usepackage[T1]{fontenc}
\usepackage[english]{babel}
\usepackage{amscd,amsmath,amsthm,amssymb,graphics}
\usepackage{lmodern,pst-node}
\usepackage{pstricks}
\usepackage{subcaption}
\usepackage{tikz}

\usetikzlibrary{positioning}
\usepackage{tikz,fullpage}
\usetikzlibrary{arrows,%
                petri,%
                topaths, automata}%
\usepackage{etex}

\usetikzlibrary{decorations.markings}
\tikzstyle{vertex}=[circle, draw, inner sep=0pt, minimum size=6pt]
\newcommand{\vertex}{\node[vertex]}

\usepackage{multicol}
\usepackage{epic,eepic}
\usepackage{amsfonts,amssymb,amscd,amsmath,enumitem,verbatim}
\psset{unit=0.7cm,linewidth=0.8pt,arrowsize=2.5pt 4}
\newpsstyle{fatline}{linewidth=1.5pt}
\newpsstyle{fyp}{fillstyle=solid,fillcolor=verylight}
\definecolor{verylight}{gray}{0.97}
\definecolor{light}{gray}{0.9}
\definecolor{medium}{gray}{0.85}
\definecolor{dark}{gray}{0.6}

\def\frk{\frak}               % font for "Fraktur"

\def\Phi{{\frk n}}
\def\Phi{{\frk N}}
\def\opn#1#2{\def#1{\operatorname{#2}}} % to make operators
\opn\chara{char} \opn\length{\ell} \opn\pd{pd} \opn\rk{rk}
\opn\projdim{proj\,dim} \opn\injdim{inj\,dim} \opn\rank{rank}
\opn\depth{depth} \opn\grade{grade} \opn\height{height}
\opn\embdim{emb\,dim} \opn\codim{codim}

\opn\Tr{Tr} \opn\bigrank{big\,rank}
\opn\superheight{superheight}\opn\lcm{lcm}
\opn\trdeg{tr\,deg}%\emph{
\opn\reg{reg} \opn\lreg{lreg} \opn\ini{in} \opn\lpd{lpd}
\opn\size{size}\opn\bigsize{bigsize}
\opn\cosize{cosize}\opn\bigcosize{bigcosize}
\opn\sdepth{sdepth}\opn\sreg{sreg}
\opn\link{link}\opn\fdepth{fdepth}
\opn\deg{deg}
\opn\max{max}
\opn\indeg{indeg}
\opn\min{min}
\opn\psln{psln}
\opn\div{div} \opn\Div{Div} \opn\cl{cl} \opn\Cl{Cl}
\let\epsilon\varepsilon
\let\phi=\varphi
\let\kappa=\varkappa

\opn\Spec{Spec} \opn\Supp{Supp} \opn\supp{supp} \opn\Sing{Sing}
\opn\Ass{Ass} \opn\Min{Min}\opn\Mon{Mon} \opn\dstab{dstab} \opn\astab{astab}
\opn\Syz{Syz}
\opn\Ann{Ann} \opn\Rad{Rad} \opn\Soc{Soc}
\opn\Im{Im}
 \opn\Ind{Ind}
 \opn\del{del}
 \opn\Ker{Ker} \opn\Coker{Coker} \opn\Am{Am}
\opn\Hom{Hom} \opn\Tor{Tor} \opn\Ext{Ext} \opn\End{End}
\opn\Aut{Aut} \opn\id{id}

\opn\nat{nat}
\opn\pff{pf}%   \pf exists already
\opn\Pf{Pf} \opn\GL{GL} \opn\SL{SL} \opn\mod{mod} \opn\ord{ord}
\opn\Gin{Gin} \opn\Hilb{Hilb}\opn\sort{sort}
\opn\initial{init}
\opn\ende{end}
\opn\height{height}
\opn\bight{bight}
\opn\hte{ht}
\opn\indeg{indeg}
\opn\reg{reg}
\opn\depth{depth}
\opn\type{type}
\opn\ldim{ldim}
\opn\maxdeg{maxdeg}
\opn\aff{aff} \opn\con{conv} \opn\relint{relint} \opn\st{st}
\opn\lk{lk} \opn\cn{cn} \opn\core{core} \opn\vol{vol}
\opn\link{link} \opn\star{star}\opn\lex{lex}
\opn\gr{gr}

\def\pot#1#2{#1[\kern-0.28ex[#2]\kern-0.28ex]}

\opn\dirlim{\underrightarrow{\lim}}
\opn\inivlim{\underleftarrow{\lim}}
\def\Implies{\ifmmode\Longrightarrow \else
        \unskip${}\Longrightarrow{}$\ignorespaces\fi}
\def\implies{\ifmmode\Rightarrow \else
        \unskip${}\Rightarrow{}$\ignorespaces\fi}
\def\iff{\ifmmode\Longleftrightarrow \else
        \unskip${}\Longleftrightarrow{}$\ignorespaces\fi}

\let\:=\colon
 \theoremstyle{plain}
\newtheorem{Theorem}{Theorem}[section]
 \newtheorem{Lemma}[Theorem]{Lemma}

 \newtheorem{Question}[Theorem]{Question}
 \theoremstyle{definition}
 \newtheorem{Definition}[Theorem]{Definition}

 \newtheorem{Example}[Theorem]{Example}

\let\epsilon\varepsilon
\let\kappa=\varkappa
\def\qed{\ifhmode\textqed\fi
      \ifmmode\ifinner\quad\qedsymbol\else\dispqed\fi\fi}
\def\textqed{\unskip\nobreak\penalty50
       \hskip2em\hbox{}\nobreak\hfil\qedsymbol
       \parfillskip=0pt \finalhyphendemerits=0}
\def\dispqed{\rlap{\qquad\qedsymbol}}

\opn\dis{dis}
\def\pnt{{\raise0.5mm\hbox{\large\bf.}}}

\opn\Lex{Lex}
\begin{document}

\author[Mafi, Naderi and Soufivand]{ Amir Mafi, Dler Naderi and Parasto Soufivand}
\title{On vertex decomposability and regularity of graphs}

\address{Amir Mafi, Department of Mathematics, University Of Kurdistan, P.O. Box: 416, Sanandaj, Iran.}
\email{A\_Mafi@ipm.ir}
\address{Dler Naderi, Department of Mathematics, University of Kurdistan, P.O. Box: 416, Sanandaj,
Iran.}
\email{dler.naderi65@gmail.com}
\address{Parasto Soufivand, Department of Mathematics, University of Kurdistan, P.O. Box: 416, Sanandaj,
Iran.}
\email{Parisoufivand@gmail.com}

\begin{abstract}
There are two motivating questions in \cite{MTS, MTS1} about Castelnuovo-Mumford regularity and vertex decomposability of simple graphs.  In this paper, we give negative answers to the questions by providing two counterexamples.
\end{abstract}


\maketitle
\section{Introduction}
Throughout this paper, we assume that $R=K[x_1,\ldots,x_n]$ is the polynomial ring over a field $K$ and suppose that $G$ is a finite simple graph on the vertex set $V=\{x_1,\ldots,x_n\}$ and the edge set $E$.
For a vertex $v$ of $G$ the set of all neighbors of $v$ is denoted by $N(v)$ and we denote by $N[v]$ the set $N(v)\cup\{v\}$ and also we denote by $\deg(v)$ the number $\mid N(v)\mid$. An independent set of $G$ is a  subset $A$ of $V(G)$ such that none of its elements are adjacent. The {\it edge ideal} of the graph $G$ is the quadratic square-free monomial ideal $I(G)=\langle x_ix_j\mid \{x_i,x_j\}\in E\rangle$ and was first introduced by Villarreal \cite{Vi}. Two edges $\{x,y\}$ and $\{z,u\}$ of $G$ are called $3$-{\it disjoint} if the induced subgraph of $G$ on $\{x,y,z,u\}$ is disconnected or equivalently in the complement of $G$ the induced graph on  $\{x,y,z,u\}$ is a four-cycle. A subset $A$ of edges of $G$ is called a pairwise $3$-disjoint set of edges in $G$ if each pair of edges of $A$ is $3$-disjoint, see \cite{Ku, MTS1, Z}. The maximum cardinality of all pairwise $3$-disjoint sets of edges in $G$ is denoted by $a(G)$, see \cite{Ku,MTS1,Z}. Note that $a(G)$ is called {\it induced matching number}. The Castelnuovo-Mumford regularity of a graded $R$-module $M$ is defined as $\reg(M)=\max\{j-i| ~\beta_{i,j}(M)\neq 0\}$. Katzmann \cite{K} proved that $\reg(R/I(G))\geq a(G)$ for every simple graph $G$.
Stanley \cite{S} defined a graded $R$-module $M$ to be {\it sequentially Cohen-Macaulay} if there exists a finite filtration of graded $R$-modules
$0=M_0\subset M_1\subset\ldots\subset M_r=M$ such that each $M_{i}/M_{i-1}$ is Cohen-Macaulay, and the Krull dimensions of the quotients are increasing:
$\dim(M_1/M_0)<\dim(M_2/M_1)<\ldots<\dim(M_r/M_{r-1})$. In particular, we call the graph $G$ sequentially Cohen-Macaulay (resp., unmixed) if $R/I(G)$ is sequentially Cohen-Macaulay (resp., unmixed).
Herzog and Hibi \cite{HH1} defined the homogeneous ideal $I$ to be {\it componentwise linear} if $(I_d)$ has a linear resolution for all $d$, where $(I_d)$ is the ideal generated by all degree $d$ elements of $I$. They proved that if $I$ is a square-free monomial ideal, then $R/I$ is sequentially Cohen-Macaulay if and only if the square-free Alexander dual $I^{\vee}$ is componentwise linear. It is known that if  $I$ has a linear resolution, then $I$ is componentwise linear. Note that for a square-free monomial ideal $I=\langle\{x_{i1}\ldots x_{in_i}\mid i=1,\ldots,t\}\rangle$ of $R$ the {\it Alexander dual} of $I$, denoted by $I^{\vee}$, is defined as $I^{\vee}=\cap_{i=1}^t\langle x_{i1},\ldots, x_{in_i}\rangle$. For a monomial ideal $I$, we write $(I_i)$ to denote the ideal generated by the degree $i$ elements of $I$.  The monomial ideal $I$ is componentwise linear if $(I_i)$ has a linear resolution for all $i$ (see \cite{HH1}). If $I$ is generated by square-free monomials, then we denote by $I_{[i]}$ the ideal generated by the square-free monomials of degree $i$ of $I$. Herzog and Hibi \cite[Proposition 1.5]{HH1} proved that the square-free monomial ideal $I$ is componentwise linear if and only if $I_{[i]}$ has a linear resolution for all $i$.

Woodroofe \cite{W} defined the graph $G$ to be {\it vertex decomposable} if it is a totally
disconnected graph (with no edges) or if the following recursive conditions hold:\\
$(i)$ there is a vertex $v$ in $G$ such that  $G\setminus v$ and $G\setminus N[v]$ are both vertex decomposable;\\
$(ii)$  no independent set in  $G\setminus N[v]$ is a maximal independent set in  $G\setminus v$.

The equality $\reg(R/I(G))=a(G)$ was proved in the following cases: $(i)$ $G$ is a tree graph; $(ii)$ $G$ is a chordal graph, where
the graph $G$ is called {\it chordal} if every cycle of length $>3$ has a chord;
$(iii)$ $G$ is a bipartite graph and unmixed;
$(iv)$ $G$ is a bipartite graph and sequentially Cohen-Macaulay;
$(v)$ $G$ is a very well-covered graph, where the graph $G$ is called
{\it very well-covered} if it is unmixed without an isolated vertices and $2\height(I(G))=\mid V\mid$;
$(vi)$ $G$ is a $C_5$-free vertex decomposable graph; $(vii)$ $G$ is an almost complete multipartite graph such that it is sequentially Cohen-Macaulay or unmixed. For details see \cite{Z, K, HV, V, MTS1, KM, H}.

Mahmoudi et al.  in \cite[Question 4.11]{MTS} and in \cite[Question 4.13]{MTS1} raised the following question:
\begin{Question}
 Let $G$ be a sequentially Cohen-Macaulay graph with $2n$ vertices which are not isolated and with $\hte (I(G)) = n$. Then do we have the following statements?
\begin{enumerate}
\item
$G$ has a vertex $v$ such that $\deg(v) = 1$.
\item
$G$ is vertex decomposable.
\item
$\reg (R/I(G)) = a(G)$.
\end{enumerate}
\end{Question}
In this paper we give a negative answer to this question by providing two counterexamples.
For every unexplained notion or terminology, we refer the reader to \cite{HH}.
%%%%%%%%%%%%%%%%%%%%%%%%%%%%%%%%%%%%%%%%%%%%%%%%%%%%%%

\section{ Counterexamples}
We start this section by recalling the following definition:
\begin{Definition}
Let $I$ be a monomial ideal of $R$ all of whose generators have degree $d$. Then $I$ has a linear resolution if for all $i\geq 0$ and for all $j\neq i+d$, $\beta_{i,j}(I)=0$. In particular, $I$ has a linear resolution if and only if $\reg(I)=d$.
\end{Definition}

\begin{Lemma}(\cite[Lemma 2.3]{AMS})\label{L}
Let $I=\langle u_1,\ldots,u_m\rangle$ be a monomial ideal with $\deg(u_i)=d_i$ and $d_i\leq d_{i+1}$ for $1\leq i\leq m-1$. If $(I_i)$ has a linear resolution for all $i<d_m$ and $\reg(I)=d_m$, then $I$ is
componentwise linear.
\end{Lemma}

By the following example we show that the Question 1.1$(1)$ and $(3)$ have  negative answers:
\begin{Example}

Let $G$ be the following graph:

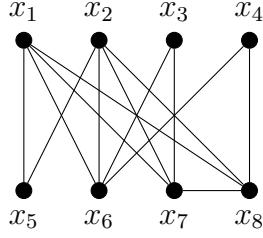
\[\begin{tikzpicture}
	\vertex[fill] (x1) at (0,2) [label=above:$x_{1}$] {};
	\vertex[fill] (x2) at (1,2) [label=above:$x_{2}$] {};
	\vertex[fill] (x3) at (2,2) [label=above:$x_{3}$] {};
\vertex[fill] (x4) at (3,2) [label=above:$x_{4}$] {};

	\vertex[fill] (x5) at (0,0) [label=below:$x_{5}$] {};
	\vertex[fill] (x6) at (1,0) [label=below:$x_{6}$] {};
	\vertex[fill] (x7) at (2,0) [label=below:$x_{7}$] {};
\vertex[fill] (x8) at (3,0) [label=below:$x_{8}$] {};

	\path
		(x1) edge (x5)
		(x1) edge (x6)
		(x1) edge (x7)
		(x1) edge (x8)
		
		(x2) edge (x5)
		(x2) edge (x6)
		(x2) edge (x7)
         (x2) edge (x8)

		(x3) edge (x6)
		(x3) edge (x7)
		
		(x4) edge (x6)
		(x4) edge (x8)
		(x7) edge (x8);
\end{tikzpicture}\]

Then we may consider the edge ideal
\[I=(x_{1}x_{5}, x_{1}x_{6}, x_{1}x_{7}, x_{1}x_{8}, x_{2}x_{5}, x_{2}x_{6}, x_{2}x_{7}, x_{2}x_{8}, x_{3}x_{6}, x_{3}x_{7}, x_{4}x_{6}, x_{4}x_{8}, x_{7}x_{8})\]
of $R=K[x_1, \ldots, x_{8}]$.
 This ideal has the following primary decomposition
\begin{align*}
I=&(x_5, x_6, x_7, x_8) \cap( x_1, x_2, x_3, x_4, x_7) \cap (x_1, x_2, x_3, x_4, x_8) \cap (x_1, x_2, x_3, x_6, x_8) \\
&\cap (x_1, x_2, x_4, x_6, x_7) \cap (x_1, x_2, x_6, x_7, x_8).
\end{align*}
So $\hte(I)=4$ and \[I^{\vee}=(x_5x_6x_7x_8,x_1x_2x_3x_4x_7,x_1x_2x_3x_4x_8,x_1x_2x_3x_6x_8,x_1x_2x_4x_6x_7,x_1x_2x_6x_7x_8).\]
Hence by using  Macaulay2 \cite{G}, we have  $\reg(R/I)=2$ and
$\reg(I^{\vee})=5$. Therefore by Lemma \ref{L} it readily follows that
 $G$ is sequentially Cohen-Macaulay.
One can easily check that for any two edges $\{x_i, x_j\}$ and $\{x_k, x_l\}$ of $G$ such that $i,j,l, k$ are different positive integers, the induced subgraph of $G$ on the vertices $\{x_i, x_j, x_k, x_l\}$ is connected. Therefore, $a(G)=1\ne\reg(R/I)$
giving a negative answer to Question 1.1.(1) and, in addition, $G$ does not have a vertex of degree $1$ contradicting Question 1.1.(3).
\end{Example}

Recall that a {\it circulant graph} is defined as follows: let $n\geq 1$ be an integer and let $S\subseteq\{1,\ldots,\lfloor\frac{n}{2}\rfloor\}$. The  circulant graph $C_n(S)$ is the graph on $n$ vertices
 $V=\{x_1,\ldots,x_n\}$ such that $\{x_i,x_j\}$ is an edge of  $C_n(S)$ if and only if $\min\{\vert i-j\vert,n-\vert i-j\vert\}\in S$. For ease of notation, we  write $C_n(a_1,\ldots,a_t)$ instead of $C_n(\{a_1,\ldots,a_t\})$, for more details see \cite{KMV}.
Let $\Delta$ be a simplicial complex on the vertex set $V = \{x_1,\ldots, x_n \}$. Members of $\Delta$ are called faces of $\Delta$ and a facet of $\Delta$ is a maximal face of $\Delta$ with respect to inclusion.  The simplicial complex $\Delta$ is pure if every facet has the same cardinality. Also, the simplicial complex $\Delta$ with the facets $F_1,\dots, F_r$ is denoted by $\Delta=\langle{F_1,\dots,F_r}\rangle$. The simplicial complex $\Delta$ is called a {\it simplex} when it has a unique facet.
For the simplicial complex $\Delta$ and the face $F \in \Delta$, one can introduce two new simplicial complexes. The {\it deletion} of $F$ from $\Delta$ is $\del_{\Delta}(F) = \{ A \in \Delta \vert F\cap A=\emptyset\} $. The {\it link} of $F$ in $\Delta$ is $\lk_{\Delta}(F) =\{ A\in  \Delta \vert F \cap A =\emptyset, A\cup F\in\Delta \}$. If $F =\{v \}$, we write $\del_{\Delta} v$ (resp. $\lk_{\Delta}v$) instead of $\del_{\Delta}(\{v\})$ (resp. $\lk_{\Delta}(\{v\})$); see \cite{HH} for details information.
The {\it Stanley-Reisner} ideal of $\Delta$ over $K$ is the ideal $I_{\Delta}$ of $R$ which is generated by those square-free monomials $x_F$ with $F \notin \Delta$, where  $x_{ F} = \prod_{x_i \in F} x_{i}$. Let $I$ be an arbitrary square-free monomial ideal.
Then there is a unique simplicial complex $\Delta$ such that $I = I_{\Delta}$.
Following \cite{W} a simplicial complex $\Delta$ is  recursively defined to be {\it vertex decomposable} if it is either a simplex or else has some vertex $v$ so that
$(i)$ both $\del_{\Delta} v$ and $\lk_{\Delta}v$ are vertex decomposable, and
$(ii)$  no face of  $\lk_{\Delta}v$ is a facet of $\del_{\Delta} v$.

A simplicial complex $\Delta$ is {\it shellable} if the facets of $\Delta$ can be ordered, say $F_1,\ldots,F_s$, such that for all $1\leq i<j\leq s$, there exists some $x\in F_j\setminus F_i$ and some $k\in\{1,2,\ldots,j-1\}$ with $F_j\setminus F_k=\{x\}$. Hence  if $\Delta$ is shellable with shelling order $F_1,\ldots,F_s$, then
for each $2\leq j\leq s$, the subcomplex $\langle F_1,\ldots,F_{j-1}\rangle\cap\langle F_j\rangle$ is pure of dimension $\dim F_j-1$, for detials see \cite[Section 8.2]{HH}. The following implications hold:\\
vertex decomposable $\Longrightarrow$ shellable $\Longrightarrow$ sequentially Cohen-Macaulay.\\
Also, both implications are known to be strict.

The independence complex of the graph $G$ is defined by $\Ind(G)=\{F\subseteq V\mid F$ is an independence set in $G\}$. It is clear $I(G)=I_{\Ind(G)}$. Let $v$ be a vertex of $G$. By \cite{H} we have the following relations:\\
$\del_{\Ind(G)} v=\Ind(G\setminus v)$ and $\lk_{\Ind(G)}v=\Ind(G\setminus N[v])$.
Therefore one can deduce that the graph $G$ is vertex decomposable if and only if the independence complex
$\Ind(G)$ is vertex decomposable.

\begin{Theorem}(\cite[Theorem 6.1 (iii)]{KMV})\label{T1}
The graph $C_{16}(1, 4, 8)$ is the smallest well-covered circulant that is shellable but not
vertex decomposable.
\end{Theorem}

By the following example we show that Question 1.1$(2)$ has a negative answer:

\begin{Example}
Let $I$ be an ideal of $R=K[x_1, \ldots, x_{26}]$ generated by the following monomials

\begin{center}
  \centering
  \scalebox{0.8}{%
\begin{tabular}{lllllllllllll}
$x_{16}x_{26}$ &   $x_{15}x_{26}$       & $x_{13}x_{26}$ & $x_{12}x_{26}$ & $x_{10}x_{26}$ &  $x_{8}x_{26}$ & $x_{7}x_{26}$ & $x_{6}x_{26}$ & $ x_{5}x_{26}$ & $x_{4}x_{26}$ & $  x_{3}x_{26}$ &$ x_{2}x_{26} $ & $ x_{1}x_{26}$\\

$x_{16}x_{25}$ &   $x_{15}x_{25}$       & $x_{13}x_{25}$ & $x_{12}x_{25}$ & $x_{10}x_{25}$ &  $x_{8}x_{25}$ & $x_{7}x_{25}$ & $x_{6}x_{25}$ & $ x_{5}x_{25}$ & $x_{4}x_{25}$ & $  x_{3}x_{25}$ &$ x_{2}x_{25} $ & $ x_{1}x_{25}$\\

$x_{16}x_{24}$ &   $x_{15}x_{24}$       & $x_{13}x_{24}$ & $x_{12}x_{24}$ & $x_{10}x_{24}$ &  $x_{8}x_{24}$ & $x_{7}x_{24}$ & $x_{6}x_{24}$ & $ x_{5}x_{24}$ & $x_{4}x_{24}$ & $  x_{3}x_{24}$ &$ x_{2}x_{24} $ & $ x_{1}x_{24}$\\

$x_{16}x_{23}$ &   $x_{15}x_{23}$       & $x_{13}x_{23}$ & $x_{12}x_{23}$ & $x_{10}x_{23}$ &  $x_{8}x_{23}$ & $x_{7}x_{23}$ & $x_{6}x_{23}$ & $ x_{5}x_{23}$ & $x_{4}x_{23}$ & $  x_{3}x_{23}$ &$ x_{2}x_{23} $ & $ x_{1}x_{23}$\\

$x_{16}x_{22}$ &   $x_{15}x_{22}$       & $x_{13}x_{22}$ & $x_{12}x_{22}$ & $x_{10}x_{22}$ &  $x_{8}x_{22}$ & $x_{7}x_{22}$ & $x_{6}x_{22}$ & $ x_{5}x_{22}$ & $x_{4}x_{22}$ & $  x_{3}x_{22}$ &$ x_{2}x_{22} $ & $ x_{1}x_{22}$\\

$x_{16}x_{21}$ &   $x_{15}x_{21}$       & $x_{13}x_{21}$ & $x_{12}x_{21}$ & $x_{10}x_{21}$ &  $x_{8}x_{21}$ & $x_{7}x_{21}$ & $x_{6}x_{21}$ & $ x_{5}x_{21}$ & $x_{4}x_{21}$ & $  x_{3}x_{21}$ &$ x_{2}x_{21} $ & $ x_{1}x_{21}$\\

$x_{16}x_{20}$ &   $x_{15}x_{20}$       & $x_{13}x_{20}$ & $x_{12}x_{20}$ & $x_{10}x_{20}$ &  $x_{8}x_{20}$ & $x_{7}x_{20}$ & $x_{6}x_{20}$ & $ x_{5}x_{20}$ & $x_{4}x_{20}$ & $  x_{3}x_{20}$ &$ x_{2}x_{20} $ & $ x_{1}x_{20}$\\

$x_{16}x_{19}$ &   $x_{15}x_{19}$       & $x_{13}x_{19}$ & $x_{12}x_{19}$ & $x_{10}x_{19}$ &  $x_{8}x_{19}$ & $x_{7}x_{19}$ & $x_{6}x_{19}$ & $ x_{5}x_{19}$ & $x_{4}x_{19}$ & $  x_{3}x_{19}$ &$ x_{2}x_{19} $ & $ x_{1}x_{19}$\\

$x_{16}x_{18}$ &   $x_{15}x_{18}$       & $x_{13}x_{18}$ & $x_{12}x_{18}$ & $x_{10}x_{18}$ &  $x_{8}x_{18}$ & $x_{7}x_{18}$ & $x_{6}x_{18}$ & $ x_{5}x_{18}$ & $x_{4}x_{18}$ & $  x_{3}x_{18}$ &$ x_{2}x_{18} $ & $ x_{1}x_{18}$\\

$x_{16}x_{17}$ &   $x_{15}x_{17}$       & $x_{13}x_{17}$ & $x_{12}x_{17}$ & $x_{10}x_{17}$ &  $x_{8}x_{17}$ & $x_{7}x_{17}$ & $x_{6}x_{17}$ & $ x_{5}x_{17}$ & $x_{4}x_{17}$ & $  x_{3}x_{17}$ &$ x_{2}x_{17} $ & $ x_{1}x_{17}$\\

$x_{15}x_{16}$ &   $x_{12}x_{16}$       & $x_{8}x_{16}$ & $x_{4}x_{16}$ & $x_{1}x_{16}$ &  $x_{14}x_{15}$ & $x_{11}x_{15}$ & $x_{7}x_{15}$ & $ x_{3}x_{15}$ & $x_{13}x_{14}$ & $  x_{10}x_{14}$ &$ x_{6}x_{14} $ & $ x_{2}x_{14}$\\

$x_{12}x_{13}$ &   $x_{9}x_{13}$       & $x_{5}x_{13}$ & $x_{1}x_{13}$ & $x_{11}x_{12}$ &  $x_{8}x_{12}$ & $x_{4}x_{12}$ & $x_{10}x_{11}$ & $ x_{7}x_{11}$ & $x_{3}x_{11}$ & $  x_{9}x_{10}$ &$ x_{6}x_{10} $ & $ x_{2}x_{10}$\\

$x_{8}x_{9}$ &   $x_{5}x_{9}$       & $x_{1}x_{9}$ & $x_{7}x_{8}$ & $x_{4}x_{8}$ &  $x_{6}x_{7}$ & $x_{3}x_{7}$ & $x_{5}x_{6}$ & $ x_{2}x_{6}$ & $x_{4}x_{5}$ & $  x_{1}x_{5}$ &$ x_{3}x_{4} $ & $ x_{2}x_{3}$\\

$x_1 x_2$
\end{tabular}}
\end{center}

 The ideal $I$ is an edge ideal of a graph, say $G$.
 This ideal has the form
 \[I=(J,x_{17}, x_{18}, \cdots, x_{26}) \cap (x_1, \cdots, x_8, x_{10}, x_{12}, x_{13}, x_{15}, x_{16}),\]
 where $J$ is the edge ideal of circulant graph $C_{16}(1, 4, 8)$.
 This ideal has the following primary decomposition
\begin{align*}
I=\mathop  \cap \limits_{i = 1}^{80} {(\frk{p}_i , x_{17}, x_{18}, \cdots, x_{26})} \cap  (x_1, \cdots, x_8, x_{10}, x_{12}, x_{13}, x_{15}, x_{16});
\end{align*}
where $\frk{p_{i}}$ for $1 \leq i \leq 80$ is an associated prime of circulant graph $C_{16}(1, 4, 8)$.
Therefore $\hte (I)=13$ and
the simplicial complex $\Ind(G)$ has $81$ facets as follows:\\

\scalebox{0.69}{%
~ $F_{0}=\{ x_{9}, x_{11}, x_{14}, x_{17},x_{18},x_{19},x_{20},x_{21},x_{22},x_{23},x_{24},x_{25},x_{26}\},$ }\\
%\end{small}
%\begin{LTR}
%\FloatBarrier
%\begin{center}
%\begin{table}[H]
%  \centering
  %\caption{}\label{tab3}
 % \vskip -.1cm
 % \end{table}
%\begin{small}
  \scalebox{0.66}{%
\begin{tabular}{lllll}

$F_{1}=\{ x_9, x_{11}, x_{14}, x_{16} \},$ &   $F_{2}=\{ x_5, x_{11}, x_{14}, x_{16} \},$       & $F_{3}=\{ x_7, x_{9}, x_{14}, x_{16} \},$ & $F_{4}=\{ x_3, x_{9}, x_{14}, x_{16} \},$ & $F_{5}=\{ x_5, x_{7}, x_{14}, x_{16} \},$ \\

$F_{6}=\{ x_3, x_{5}, x_{14}, x_{16} \},$ &   $F_{7}=\{ x_6, x_{9}, x_{11}, x_{16} \},$       & $F_{8}=\{ x_5, x_{7}, x_{10}, x_{16} \},$ & $F_{9}=\{ x_2, x_{5}, x_{11}, x_{16} \},$ & $F_{10}=\{ x_2, x_{9}, x_{11}, x_{16} \},$ \\

$F_{11}=\{ x_2, x_{7}, x_{13}, x_{16} \},$ &   $F_{12}=\{ x_7, x_{10}, x_{13}, x_{16} \},$       & $F_{13}=\{ x_2, x_{11}, x_{13}, x_{16} \},$ & $F_{14}=\{ x_6, x_{11}, x_{13}, x_{16} \},$ & $F_{15}=\{ x_3, x_{5}, x_{10}, x_{16} \},$ \\
$F_{16}=\{ x_3, x_{10}, x_{13}, x_{16} \},$ &   $F_{17}=\{ x_3, x_{6}, x_{13}, x_{16} \},$       & $F_{18}=\{ x_2, x_{7}, x_{9}, x_{16} \},$ & $F_{19}=\{ x_3, x_{6}, x_{9}, x_{16} \},$ & $F_{20}=\{ x_2, x_{5}, x_{7}, x_{16} \},$ \\

$F_{21}=\{ x_7, x_{9}, x_{12}, x_{14} \},$ &   $F_{22}=\{ x_1, x_{4}, x_{10}, x_{15} \},$       & $F_{23}=\{ x_1, x_{8}, x_{10}, x_{15} \},$ & $F_{24}=\{ x_5, x_{8}, x_{10}, x_{15} \},$ & $F_{25}=\{ x_1, x_{10}, x_{12}, x_{15} \},$ \\

$F_{26}=\{ x_4, x_{10}, x_{13}, x_{15} \},$ &   $F_{27}=\{ x_8, x_{10}, x_{13}, x_{15} \},$       & $F_{28}=\{ x_5, x_{10}, x_{12}, x_{15} \},$ & $F_{29}=\{ x_3, x_{9}, x_{12}, x_{14} \},$ & $F_{30}=\{ x_3, x_{8}, x_{10}, x_{13} \},$ \\
$F_{31}=\{ x_3, x_{5}, x_{8}, x_{14} \},$ &   $F_{32}=\{ x_5, x_{8}, x_{11}, x_{14} \},$       & $F_{33}=\{ x_6, x_{8}, x_{11}, x_{13} \},$ & $F_{34}=\{ x_6, x_{8}, x_{13}, x_{15} \},$ & $F_{35}=\{ x_4, x_{6}, x_{13}, x_{15} \},$ \\

$F_{36}=\{x_2, x_{8}, x_{13}, x_{15} \},$ &   $F_{37}=\{x_2, x_{8}, x_{11}, x_{13} \},$       & $F_{38}=\{x_1, x_{4}, x_{6}, x_{15} \},$ & $F_{39}=\{  x_4, x_{6}, x_{9}, x_{15} \},$ & $F_{40}=\{ x_6, x_{9}, x_{12}, x_{15} \},$ \\

$F_{41}=\{  x_1, x_{6}, x_{12}, x_{15} \},$ &   $F_{42}=\{ x_1, x_{6}, x_{8}, x_{15}  \},$       & $F_{43}=\{ x_2, x_{4}, x_{13}, x_{15} \},$ & $F_{44}=\{ x_2, x_{9}, x_{12}, x_{15} \},$ & $F_{45}=\{ x_2, x_{4}, x_{9}, x_{15} \},$ \\

$F_{46}=\{ x_4, x_{6}, x_{11}, x_{13} \},$ &   $F_{47}=\{ x_4, x_{9}, x_{11}, x_{14} \},$       & $F_{48}=\{ x_4, x_{7}, x_{9}, x_{14} \},$ & $F_{49}=\{ x_2, x_{4}, x_{11}, x_{13} \},$ & $F_{50}=\{x_5, x_{7}, x_{10}, x_{12}  \},$ \\

$F_{51}=\{x_1, x_{3}, x_{8}, x_{14}   \},$ &   $F_{52}=\{x_1, x_{8}, x_{11}, x_{14}   \},$       & $F_{53}=\{ x_1, x_{3}, x_{12}, x_{14} \},$ & $F_{54}=\{ x_1, x_{7}, x_{12}, x_{14} \},$ & $F_{55}=\{ x_1, x_{7}, x_{10}, x_{12}\},$ \\

$F_{56}=\{ x_3, x_{6}, x_{8}, x_{13}  \},$ &   $F_{57}=\{ x_5, x_{7}, x_{12}, x_{14} \},$       & $F_{58}=\{ x_3, x_{5}, x_{12}, x_{14} \},$ & $F_{59}=\{ x_3, x_{5}, x_{10}, x_{12} \},$ & $F_{60}=\{ x_1, x_{3}, x_{10}, x_{12} \},$ \\

$F_{61}=\{x_2,  x_7,  x_9,  x_{12} \},$ &   $F_{62}=\{ x_3, x_{6}, x_{9}, x_{12}   \},$       & $F_{63}=\{x_2, x_{5}, x_{7}, x_{12}   \},$ & $F_{64}=\{ x_2, x_{5}, x_{8}, x_{11}\},$ & $F_{65}=\{x_1, x_{6}, x_{8}, x_{11} \},$ \\

$F_{66}=\{x_2, x_{4}, x_{9}, x_{11} \},$ &   $F_{67}=\{ x_4, x_{6}, x_{9}, x_{11} \},$       & $F_{68}=\{   x_1, x_{3}, x_{6}, x_{12}\},$ & $F_{69}=\{ x_2, x_{5}, x_{8}, x_{15} \},$ & $F_{70}=\{x_2, x_{5}, x_{12}, x_{15}     \},$ \\

$F_{71}=\{ x_1, x_{4}, x_{6}, x_{11} \},$ &   $F_{72}=\{ x_1, x_{4}, x_{11}, x_{14} \},$ & $F_{73}=\{ x_1, x_{4}, x_{7}, x_{14} \},$ & $F_{74}=\{ x_3, x_{5}, x_{8}, x_{10} \},$ & $F_{75}=\{ x_1, x_{3}, x_{8}, x_{10} \},$ \\

$F_{76}=\{ x_1, x_{4}, x_{7}, x_{10} \},$ &   $F_{77}=\{ x_4, x_{7}, x_{10}, x_{13} \},$ & $F_{78}=\{ x_1, x_{3}, x_{6}, x_{8} \},$ & $F_{79}=\{ x_2, x_{4}, x_{7}, x_{9} \},$ & $F_{80}=\{ x_2, x_{4}, x_{7}, x_{13} \}$ \\
\end{tabular}}
%\end{small}
%\end{center}
%\end{LTR}
%\vspace{5}
\\

By the proof of Theorem \ref{T1}, we have $F_{1}, \ldots, F_{80} $ is a shelling order of $\Ind(C_{16}(1, 4, 8))$ and the graph $C_{16}(1, 4, 8)$ is the smallest well-covered circulant that is shellable but not
vertex decomposable. We claim that $F_{0}, F_{1}, \ldots, F_{80}$ is a shelling order of $\Ind(G)$.
Since $F_{1}, \ldots, F_{80} $ is a shelling order, it is enough to show that for each $i$, there exists some $v \in F_{i} \setminus F_{0}$ and some $k < i$ such that $F_{i} \setminus F_{k} =\{ v\}$. If  $i=1$, then it is clear
$F_{1} \setminus F_{0} =\{ x_{16}\}$. Now we assume that $1\ne i\leq 80$. Since $F_{i} \setminus F_{1}\subseteq
F_{i} \setminus F_{0}$, we may choose $v \in F_{i} \setminus F_{1}$ and so there exists some $1\leq k < i$ such that $F_{i} \setminus F_{k} =\{ v\}$. Therefore $\Ind(G)$ is shellable and so $G$ is sequentially Cohen-Macaulay.

Now, we claim that for each element $x_t $ with $1\leq t\leq 26$, $\del_{\Ind(G)} (x_t)$ is not vertex decomposable.
If  $x_{t} \in  \{x_9, x_{11}, x_{14}, x_{17},\ldots,x_{26} \}$, then by using the definition on the above facets it is obvious that  $\del_{\Ind(G)} (x_t)$ has a facet, say $F^{\prime}$, such that $F^{\prime} \ne F_{i}$ for $0 \leq i \leq 80$, and in this case $\del_{\Ind(G)} (x_t)$  is not vertex decomposable. For the remaining claim, we assume that  $x_{t} \in \{ x_{1}, \ldots, x_{8}, x_{10}, x_{12}, x_{13}, x_{15}, x_{16}\}$ and we will show that $\del_{\Ind(G)} (x_t)$ is not shellable and so it is not vertex decomposable.
By contrary, let  $\del_{\Ind(G)} (x_t)$ be shellable and so we may consider the shelling order  $F_0=F_{s_0}, F_{s_1}, \ldots, F_{s_r}$. By this shelling order we have $F_{0}=(F_{s_1}\setminus  \{x_{m} \})\cup \{ x_{17},\ldots,x_{26}\}$ for some $x_m \in F_{s _1}$  and for all $i$ and $j <i$ there exists $x_{l} \in F_{s_i} \setminus F_{s_j}$ and $k < i$ such that $F_{s_i} \setminus F_{s_k} =\{ x_{l} \}$. By this assumption we claim that
$F_{s_1}, \ldots, F_{s_r}$ is  shellable and for this it is enough for such $k$ to assume
$F_{s_k}=F_{0}$. In this case  $F_{s_i}=(F_{0}\setminus \{ x_{17},\ldots,x_{26}\})\cup  \{x_{l} \}=\{ x_{9}, x_{11}, x_{14}, x_{l} \}$. We may assume  $F_{s_i}\ne F_{s_1}$.  Since $F_{s_i}=\{ x_{9}, x_{11}, x_{14}, x_{l} \}$
and $F_{s_1}=\{x_9,x_{11},x_{14},x_m\}$, we have  $F_{s_i}\setminus F_{s_1}=\{x_l\}$.
It therefore follows that $ F_{s_1}, \ldots, F_{s_r}$ is a shelling order. Hence  $\del_{\Ind(C_{16}(1, 4, 8))} (x_t)=\langle F_{s_1}, \ldots, F_{s_r} \rangle$ and this means that $\del_{\Ind(C_{16}(1, 4, 8))} (x_t)$ is pure  shellable  and Cohen-Macaulay.
This is a contradiction by the proof of Theorem \ref{T1}. Thus $\del_{\Ind(G)} (x_t)$ is not shellable  and so $G$ is not vertex decomposable. Hence we construct a sequentially Cohen-Macaulay graph with $26$ vertices such that $\hte(I)=13$ but it is not vertex decomposable.
\end{Example}

{\bf Acknowledgments:} The authors are indebted to Adam Van Tuyl for suggestion and many valuable comments.
We also thank the referee for a careful reading of the paper and for the improvements suggested.

%%%%%%%%%%%%%%%%%%%%%%%%%%%%%%%%%%%%%%%%%%%%%%%%%%%%%%%%%%%%%%%%%%%%%%%%%%%

\end{document}